\documentclass[12pt,draft]{article}

\usepackage{amsfonts,amsmath,amssymb,theorem}

\usepackage[english]{babel}

\DeclareMathOperator{\loc}{loc}

\DeclareMathOperator{\esssup}{esssup}

\begin{document}

\centerline{\bf BRENNAN'S CONJECTURE FOR COMPOSITION OPERATORS}
\centerline{\bf ON SOBOLEV SPASES}

\vskip 0.3cm
\centerline{\bf V.~Gol'dshtein and A.~Ukhlov}
\vskip 0.3cm

\vskip 0.3cm

\vskip 0.3cm

\noindent
{\bf Abstract.}
{We show that Brennan's conjecture is equivalent to boundedness of composition operators on homogeneous Sobolev spaces, that are generated by conformal homeomorphisms of simply connected plane domains to the unit disc. A geometrical interpretation of Brennan's conjecture in terms of integrability of $p$-distortion is given.  }

\vskip 0.5cm

{\bf 1 \,Introduction}
\vskip 0.3cm

The purpose of this paper is to point out a connection between Brennan's conjecture about integrability of derivatives of plane conformal mappings and boundedness of composition operators of homogeneous Sobolev spaces $L^1_p(\mathbb D)$ to $L^1_q(\Omega)$ generated by conformal homeomorphisms $\varphi : \Omega\to\mathbb D$ of a simply connected plane domain with nonempty boundary $\Omega\subset\mathbb R^2$ to the unit disc $\mathbb D\subset \mathbb R^2$. Brennan's conjecture goes back  to the approximation problem for analytic functions (see, for example, \cite{Br,Kl,MKh}) and is under intensive study at the present time. 

We prove that {\it any conformal homeomorphism $\varphi : \Omega\to\mathbb D$ generates by the composition rule $\varphi(f)=f\circ\varphi$ a bounded composition operator
$$
\varphi^{\ast}: L^1_{p}(\mathbb D)\to L^1_{q(p,s)}(\Omega)
$$
for any $p\in (2;+\infty)$ and $q(p,s)=ps/(p+s-2)$.} Here the number $s$ is the summability exponent for Brennan's conjecture.

This connection of Brennan's conjecture with Sobolev spaces is new and possible applications to Sobolev type embedding theorems and elliptic boundary value problems are under investigation. 

\vskip 0.3cm

{\bf Brennan's conjecture. } Brennan's conjecture concerns integrability of derivatives of plane conformal homeomorphisms $\varphi: \Omega  \to \mathbb D$ of a simply connected plane domain with nonempty boundary $\Omega\subset\mathbb R^2$ to the unit disc $\mathbb D\subset\mathbb R^2$. 

The conjecture \cite{Br} is that
\begin{equation}
\int\limits_{\Omega}|\varphi'(z)|^s~d\mu<+\infty,\quad \text{for all}\quad \frac{4}{3}<s<4.
\end{equation}
 For ${4/}{3} < s < 3$ it is a comparatively easy consequence of the Koebe distortion theorem (see, for example, \cite{Br}).

 J.~Brennan \cite{Br} (1973) extended this range to $4/3 < s < 3+\delta$ where $\delta > 0$,
and conjectured it to hold for $4/3 < s < 4$. The example $\Omega = \mathbb R^2 \setminus (-\infty,-1/4]$
shows that this range of $s$ cannot be extended. The upper bound of those $s$ for
which (1.1) is known to hold has been increased by Ch. Pommerenke to $s\leq 3.399$ by
D. Bertilsso to $s\leq 3.421$ and then to $p=3.752$ by Hedenmalm and Shimorin (2005). These results and an additional information can
be found in \cite{Ber,Pom}.

The Brennan's conjecture is correct for some special classes of domains: starlike domains, bounded domains which boundaries that are locally graphs of continuous functions etc..

\vskip 0.3cm
{\bf 2 \,Composition operators on homogeneous Sobolev spaces}
\vskip 0.3cm

 Let $\Omega\subset\mathbb R^n$, $n\geq 2$, be an Euclidean domain.
The homogeneous Sobolev space
$L^1_p(\Omega)$, $1\leq p\leq\infty$, is the space of all locally summable, weakly differentiable functions $f:\Omega\to\mathbb R$ with the finite seminorm
$$
\|f\mid L^1_p(\Omega)\|=\biggl(\int\limits_{\Omega}|\nabla f(z)|^p~d\mu\biggr)^{\frac{1}{p}},\,\,\,1\leq p<\infty,\,\,\,
\|f\mid L^1_{\infty}(\Omega)\|=\esssup\limits_{z\in\Omega}|\nabla f(z)|.
$$

Let $\varphi : \Omega\to\Omega'$ be a homeomorphism of Euclidean domains $\Omega,\Omega'\subset\mathbb R^n$. We say that $\varphi$ generates a bounded composition operator
 $\varphi^{\ast}: L^1_p(\Omega')\to L^1_q(\Omega)$, $1\leq q\leq p\leq\infty$ by the composition rule $\varphi^{\ast}(f)=f\circ\varphi$, if 
for any $f\in L^1_p(\Omega')$ the function $\varphi^{\ast}(f)\in L^1_q(\Omega)$ and there exists a constant $K<\infty$ such that
$$
\|\varphi^{\ast}(f)\mid L^1_q(\Omega)\|\leq K \|f\mid L^1_p(\Omega')\|.
$$

\vskip 0.3cm

{\bf Composition problem for homogeneous Sobolev spaces}
{\it is the problem about description of homeomorphisms $\varphi: \Omega\to\Omega'$ that generate bounded composition operators from $L^1_p(\Omega')$ to $L^1_q(\Omega)$, $1\leq q\leq p\leq\infty$.}

The quasiconformal methods to study of the composition problem arises in the  Reshetnyak's problem stated in 1968. It is the problem about description of  all isomorphisms of homogeneous Sobolev spaces $L^1_n(\Omega')$ and $L^1_n(\Omega)$ which are generated by quasiconformal mappings of Euclidean domains.  This problem was solved in \cite{VG1}, where  it was proved that a mapping $\varphi: \Omega\to \Omega'$ that maps a domain $\Omega \subset \mathbb R^n$ to a domain  $\Omega'\subset \mathbb R^n$ generates by the composition rule $\varphi^{\ast}(f) = f \circ\varphi$ an isomorphism of $L^1_n(\Omega)$ and $L^1_n(\Omega')$ if and only if $\varphi$ is a quasiconformal homeomorphism. 

The homeomorphisms that generate bounded composition operators of Sobolev spaces $L^1_1(\Omega')$ and $L^1_1(\Omega)$ were introduced by V.~G. Maz'ya in \cite{M1} as a class of
sub-areal mappings.  This pioneering work
established a connection between geometrical properties
of homeomorphisms and corresponding Sobolev spaces.
This study was extended in the book
\cite{MSh}, where the theory of multipliers was
applied to the change of variable problem in Sobolev spaces. In the framework of the quasiconformal approach to Composition problem it was proved in  \cite{VG2} that homeomorphisms that generate isomorphisms of Sobolev spaces $W^1_p(\Omega')$ and $W^1_p(\Omega)$, $p>n$, are quasiisometric. The composition problem was studied for Sobolev spaces $W^1_p(\Omega')$ and $W^1_p(\Omega)$, $n-1<p<n$, in \cite{GRo}, and Sobolev spaces $W^1_p(\Omega')$ and $W^1_p(\Omega)$, $1\leq q<n$, in \cite{Mar}. Analytic characteristics of homeomorphisms that generate bounded composition operators $\varphi^{\ast}: L^1_p(\Omega')\to L^1_p(\Omega)$, $1\leq p\leq\infty$, were studied in \cite{V1}.
The geometric description of homeomorphisms that generate bounded composition operators $\varphi^{\ast}: L^1_p(\Omega')\to L^1_p(\Omega)$, $n-1<p<\infty$ was obtained in \cite{GGR}.

New problems arise when we study composition operators on Sobolev spaces with decreasing integrability of first weak derivatives.  This problem was studied firstly in \cite{U1}. Note, that in this case a significant role in the description of the composition operators play (quasi)additive set functions, defined on open sets.  This type of composition operators has applications in the Sobolev embedding theory \cite{GGu,GU}.

The main result of \cite{U1} is

{\bf Theorem A.} {\it A homeomorphism $\varphi : \Omega\to\Omega'$, $\Omega, \Omega'\subset\mathbb R^n$, generates a bounded composition operator
$$
\varphi^{\ast} : L^1_p(\Omega')\to L^1_q(\Omega),\,\,\,1\leq q< p<\infty,
$$
if and only if $\varphi\in W^1_{1,\loc}(\Omega)$, has finite distortion, and
$$
K_{p,q}(\varphi;\Omega)=
\biggl(\int\limits_{\Omega}\biggl(\frac{|D\varphi(z)|^p}{|J(z,\varphi)|}\biggr)^{\frac{q}{p-q}}~d\mu\biggr)^{\frac{p-q}{pq}}<\infty.
$$
}

Here $D\varphi(z)$ is a formal Jacobi matrix of $\varphi$ at $z\in\Omega$ and $J(z,\varphi)=\det D\varphi(z)$ its Jacobian.
The norm $|D\varphi(z)|$ of the matrix is the norm of the linear operator defined
by this matrix in the Euclidean space $\mathbb R^n$.
In \cite{GGu} this class of mappings was studied in a connection with the Sobolev type embedding theorems. A detailed  study of composition operators on Sobolev spaces was given in \cite{VU1}.
The composition operators on Sobolev spaces in the limit case $p=\infty$ were studied in \cite{GU1, GU2}.

If $1\leq q<p=n$ this integral condition of Theorem A can be expressed in terms of integrability of the conformal distortion of the mapping $\varphi$:
$$
K(x,\varphi)=\frac{|D\varphi(z)|^n}{|J(z,\varphi)|}\,\,\,\text{a.~e. in}\,\,\,\Omega.
$$
Existence of conformal distortion almost everywhere in $\Omega$ means that $f$ is a mapping of finite distortion i.~e.  $|D\varphi(z)|=0$ for almost all $z\in Z=\{z\in\Omega: J(z,\varphi)=0\}$.

The theory of mappings of finite distortion is under intensive development at the last decades. In series of works geometrical and topological properties of these mappings were studied (see, e.g. \cite{HKos, ISv, KKMa1, KKMa2,  KKMa3, MRSY}). 

The generalization of conformal distortion is the notion of so-called $p$-distortion
$$
K_p(z,\varphi)=\frac{|D\varphi(z)|^p}{|J(z,\varphi)|}\,\,\,\text{a.~e. in}\,\,\,\Omega.
$$
introduced in \cite{GGR}. On this way mappings with a finite characteristics $K_{p,q}(\varphi;\Omega)$ are natural generalizations of quasiconformal mappings and we call them as $(p,q)$-quasiconformal homeomorphisms \cite{VU} or mappings with bounded $(p,q)$-distortion \cite{UV1}.

Composition operators on weighted Sobolev spaces were studied in \cite{UV2}. These operators have applications in weighted Sobolev type embeddings and for degenerate elliptic equations \cite{GU, GMU}.

In the recent decade an interest to composition operators on Sobolev spaces increased \cite{HenK,Klp}. For example, in \cite{Klp} was proved (by a different method than in \cite{U1}) that $(p,q)$-quasiconformal homeomorphisms generate a bounded composition operators on Sobolev spaces. 

\vskip 0.3cm

{\bf 3 \,Composition theorem for conformal homeomorphisms}

\vskip 0.3cm

Let us rewrite Theorem A in the particular case of conformal homeomorphisms of plane domains $\Omega$ and $\Omega'$. In this case $|D\varphi(z)|$ is equal to the absolute value of complex derivative $|\varphi'(z)|$ and $|\varphi'(z)|^2=J(z,\varphi)>0$ for any conformal homeomorphisms $w=\varphi(z) : \Omega\to\Omega'$. Note, that conformal homeomorphisms are smooth mappings, hence they belong to $W^1_{1,\loc}(\Omega)$ and have finite distortion.

{\bf Composition Theorem.} {\it A conformal homeomorphism $\varphi: \Omega\to\Omega'$ of plane domains $\Omega, \Omega'\subset\mathbb R^2$ generates a bounded composition operator
$$
\varphi^{\ast}: L^1_{p}(\Omega')\to L^1_q(\Omega), \,\,\,1\leq q<p\leq\infty,
$$
if and only if
$$
\int\limits_{\Omega}|\varphi'(z)|^{\frac{(p-2)q}{p-q}}~d\mu<\infty \,\,
\biggl(\int\limits_{\Omega}|\varphi'(z)|^q~d\mu<+\infty \,\, \text{for}\,\, p=\infty\biggr).
$$}
\vskip 0.3cm

{\bf Remarks:} 

{\bf 1.} The theorem was proved for homeomorphisms in the case  $1\leq q< p < \infty$ in 
\cite{U1} and for the limit case $p=\infty$ in  \cite{GU1,GU2}. 

{\bf 2.} The case $p=q$ does not interesting for the present study because absence of connections with the Brennan's conjecture. 

{\bf 3.}  Recall, that a conformal homeomorphism $\varphi: \Omega \to \Omega'$ generates an isometry $\varphi^{\ast} :  L^1_2(\Omega')\to L^1_2(\Omega)$ because $|\varphi'(z)|^2=J(z,\varphi)$. 

Indeed, in this case by the direct calculation for $f\in L^1_2(\Omega')$ we have:
\begin{multline}
\|\varphi^{\ast}f\mid L^1_2(\Omega)\|=\biggl(\int\limits_{\Omega}|\nabla (f\circ\varphi(z))|^2~d\mu\biggr)^{\frac{1}{2}}= 
\biggl(\int\limits_{\Omega}|\nabla f|^2(\varphi(z))\cdot |\varphi'(z)|^2~d\mu\biggr)^{\frac{1}{2}}\\
=\biggl(\int\limits_{\Omega}|\nabla f|^2(\varphi(z))\cdot J(z,\varphi)~d\mu\biggr)^{\frac{1}{2}}
=\biggl(\int\limits_{\Omega'}|\nabla f|^2(w)~dw\biggr)^{\frac{1}{2}}=\|f\mid L^1_2(\Omega')\|.
\nonumber
\end{multline}

\vskip 0.3cm

{\bf Equivalence of Brennan's conjecture and Composition problem for conformal homeomorphisms.} 
The Riemann Mapping Theorem states existence of a conformal homeomorphism $\varphi$ of a simply connected plane domain $\Omega$ with non\-em\-pty boundary  onto the unit disc $\mathbb D\subset \mathbb R^2$. Therefore the Composition problem for conformal homeomorphisms $\varphi: \Omega \to \Omega'$  of simply connected plane domains with non\-em\-pty boundary can be reduced to case $\Omega'=\mathbb D$.

We prove that Brennan's conjecture is equivalent to boundedness of composition operators generated by conformal homeomorphisms $\varphi: \Omega\to\mathbb D$.

\vskip 0.3cm

{\bf Equivalence Theorem.} 
{\it Brennan's conjecture (1) holds for a number $s\in ({4}/{3};4)$ if and only if any conformal homeomorphism $\varphi : \Omega\to\mathbb D$ generates a bounded composition operator
$$
\varphi^{\ast}: L^1_{p}(\mathbb D)\to L^1_{q(p,s)}(\Omega)
$$
for any $p\in (2;+\infty)$ and $q(p,s)=ps/(p+s-2)$.}
\vskip 0.3cm

{\sc Proof.} By Composition Theorem
$$
\int\limits_{\Omega}|\varphi'(z)|^{\frac{(p-2)q}{p-q}}~d\mu<\infty
$$
if and only if $\varphi: \Omega\to \mathbb D$ generates a bounded composition operator
$$
\varphi^{\ast}: L^1_{p}(\mathbb D)\to L^1_q(\Omega),\quad 1\leq q< p\leq\infty.
$$
If a number $s\in ({4}/{3};4)$ satisfies to Brennan's conjecture then  
$$
\int\limits_{\Omega}|\varphi'(z)|^s~d\mu<\infty.
$$ 
From the other side 
$$
\int\limits_{\Omega}|\varphi'(z)|^s~d\mu=\int\limits_{\Omega}|\varphi'(z)|^{\frac{(p-2)q}{p-q}}~d\mu
$$
for ${(p-2)q}/{(p-q)}=s$ i.~e. for $q:=q(p,s)=ps/(p+s-2)$. Let us check the condition of Composition Theorem that $q(p,s)=ps/(p+s-2)<p$. Because $p>2 $ we have $p+s-2>s>0$. Hence $s/(p+s-2)<1$ and $q(p,s)<p$.
The theorem proved.
\vskip 0.3cm

{\bf Remarks:} 

{\bf 1.} The restriction $p\in (2;+\infty)$ is sharp. If  $p<2$ then $q(p,s)>p$ and we have no composition operators. If $p=2$ then $q(p,s)=2$ and $\varphi'\in L^2(\Omega)$ if and only if $|\Omega|<\infty$.

 {\bf 2.} The result is correct also for $p=\infty$ and $q=s$.

\vskip 0.3cm

{\bf 4 \,Geometric versions Brennan's conjecture }

\vskip 0.3cm

{\bf Brennan's conjecture for $p$-distortion.} Define the $p$-distortion of conformal homeomorphisms as 
$$
K_p(z,\varphi)=\frac{|\varphi'(z)|^p}{J(z,\varphi)}=|\varphi'(z)|^{p-2}\,\,\, \text{for any}\,\,\, p\in[1,\infty).
$$

Using the notion of $p$-distortion and motivated by the Composition Theorem and the  Equivalence Theorem we shall reformulate the Brennan's conjecture in the following equivalent form:

\vskip 0.3cm

{\bf Geometric Brennan's conjecture.} The conjecture concerns  integrability of $p$-di\-la\-ta\-tion of plane conformal homeomorphisms $\varphi: \Omega  \to \mathbb D$ of a simply connected plane domain with nonempty boundary $\Omega\subset\mathbb R^2$ to the unit disc $\mathbb D\subset\mathbb R^2$.

The conjecture is that the $p$-distortion is integrable in degree $\alpha$ for all  
$  \frac{4}{3(p-2)}<\alpha<\frac{4}{p-2}$ if $2<p<+\infty$ and
$  \frac{-4}{2-p}<\alpha<\frac{-4}{3(2-p)}$ if $1\leq p<2$.

\vskip 0.3cm

In the context of Composition problem this reformulation of Brennan's conjecture will be useful:

{\bf Inverse Brennan's conjecture}
Brennan's conjecture can be reformulated also for conformal homeomorphisms $\psi=\varphi^{-1}:\mathbb D\to\Omega$.

The conjecture becomes
$$
\int\limits_{\mathbb D}|\psi'(w)|^r~dw<+\infty,\quad \text{for all}\quad -2<r<\frac{2}{3}.
$$
Inverse Brennan's conjecture is proved for the range $-1.78 \leq r < 2/3$ (S. Shimorin (2005)).

The next theorem demonstrates connection between Inverse Brennan's conjecture and Composition Problem  for conformal maps $\psi:D\to \Omega$.
\vskip 0.3cm

{\bf Inverse Composition Theorem.} {\it
A diffeomorphism $w=\varphi: \Omega\to\Omega'$ of plane domains $\Omega,\Omega'\subset\mathbb R^2$ generates a bounded composition operator
$$
\varphi^{\ast}: L^1_{p}(\Omega')\to L^1_q(\Omega), \,\,\, 1< q< p<\infty,
$$
if and only if  $\varphi^{-1}: \Omega'\to\Omega$ generates a bounded composition operator 
$$
(\varphi^{-1})^{\ast}: L^1_{q'}(\Omega)\to L^1_{p'}(\Omega')
$$
where ${1}/{p}+{1}/{p'}=1,\quad {1}/{q}+{1}/{q'}=1.$}

\vskip 0.3cm

{\sc Proof.} 
Let $\varphi$ generates the composition operator $\varphi^{\ast}: L^1_{p}(\Omega')\to L^1_q(\Omega)$. By Composition Theorem 
$$\int\limits_{\Omega}\biggl(\frac{|\varphi'(z)|^{p}}{|J(z,\varphi)|}\biggr)^{\frac{q}{p-q}}~d\mu<\infty.$$ Using
the equality 
$
|(\varphi^{-1})'(\varphi(z))|=|\varphi'(z)|/|J(z,\varphi)|
$
and the change of variable formula we obtain
\begin{multline}
\int\limits_{\Omega'}\biggl(\frac{|(\varphi^{-1})'(w)|^{q'}}{|J(w,\varphi^{-1})|}\biggr)^{\frac{p'}{q'-p'}}~dw=
\int\limits_{\Omega'}\biggl(\biggl(\frac{|\varphi'(\varphi^{-1}(w))|}{|J(\varphi^{-1}(w),\varphi)|}\biggr)^{q'}
\frac{1}{|J(w,\varphi^{-1})|}\biggr)^{\frac{p'}{q'-p'}}~dw\\
=\int\limits_{\Omega'}\biggl(\frac{|\varphi'(\varphi^{-1}(w))|^{q'}}{|J(\varphi^{-1}(w),\varphi)|^{q'-1}}
\biggr)^{\frac{p'}{q'-p'}}dw=
\int\limits_{\Omega'}\frac{|\varphi'(\varphi^{-1}(w))|^{\frac{pq}{p-q}}}{|J(\varphi^{-1}(w),\varphi)|^{\frac{p}{p-q}}}~dw\\
=\int\limits_{\Omega}\frac{|\varphi'(z)|^{\frac{pq}{p-q}}}{|J(z,\varphi)|^{\frac{p}{p-q}}}|J(z,\varphi)|~d\mu=\int\limits_{\Omega}\biggl(\frac{|\varphi'(z)|^{p}}{|J(z,\varphi)|}\biggr)^{\frac{q}{p-q}}~d\mu<\infty.
\nonumber
\end{multline}
Hence, by Composition Theorem the composition operator 
$$
(\varphi^{-1})^{\ast}: L^1_{q'}(\Omega)\to L^1_{p'}(\Omega'),\quad {1}/{p}+{1}/{p'}=1,\quad {1}/{q}+{1}/{q'}=1.
$$
is bounded.
The inverse assertion can be proved in the same way. The theorem proved.
\vskip 0.3cm

{\bf Remarks:} 

{\bf 1.} The result is correct also for $p=\infty$ and $q=1$.

{\bf 2.}  By the Inverse Composition Theorem the Inverse Brennan's conjecture (for conformal homeomorphisms $ \varphi: \mathbb D \to \Omega $, see, for example, \cite{Ber}) is equivalent to the Inverse Composition Problem for $p\in (1,2]$.

\vskip 0.3cm

{\bf Brennan's conjecture for mappings with bounded $(p,q)$-distortion.}

Brennan's conjecture means than conformal homeomorphisms are mappings with bo\-und\-ed $(p,q)$-distortion ($(p,q)$-quasiconformal homeomorphisms) for 
any $p\in (2;+\infty)$ and $q(p,s)=ps/(p+s-2)$. This fact gives us examples of mappings with bounded $(p,q)$-distortion.
The Brennan's conjecture in geometric terms of integrability of the $p$-distortion looks more natural for the geometric function theory and applications to PDE.

Let us recall finally that Equivalence theorem is proved only for $2<p<\infty$  and for $1<p<2$ we have no information.

\vskip 0.5cm

\noindent
Vladimir Gol'dshtein  \,\,\,\,  \hskip 4cm Alexander Ukhlov

\noindent
Department of Mathematics   \hskip 3.1cm Department of Mathematics

\noindent
Ben-Gurion University of the Negev   \hskip 1.7cm Ben-Gurion University of the Negev

\noindent
P.O.Box 653, Beer Sheva, 84105, Israel.  \hskip 1.0cm P.O.Box 653, Beer Sheva, 84105, Israel.

\noindent
E-mail: vladimir@bgu.ac.il  \hskip 3.35cm E-mail: ukhlov@math.bgu.ac.il

\end{document}